\documentclass[12pt,leqno]{amsart}
\usepackage{amsmath,amsthm}
\usepackage{amsfonts}
 \newtheorem{res}{Result}[section]

 \newtheorem{theorem}[res]{Theorem}
\newtheorem{remark}[res]{Remark}
 
 \newtheorem{lemma}[res]{Lemma}
 \newtheorem{corollary}[res]{Corollary}
 \newtheorem{definition}[res]{Definition}
\newtheorem{conjecture}[res]{Problem}
\newtheorem{example}[res]{Example}
\numberwithin{equation}{section}

\def\al{\alpha}
\def\tLam{\tilde \Lambda}
\def\Lam{\Lambda}
\def\M{{\mathcal M}}
\def\Mprod{{\mathbb M}}

\def\bbQ{\mathbb Q}

\def\F{{\mathcal F}}

\def\bid{\pi}
\def\tbid{\tilde\bid}
\def\tz{\tilde{z}}

\def\solvD{{K}}
\def\tradeD{{-K}}
\def\trade{{-{\mathcal K}}}
\def\K{\trade}
\def\hatK{{-\hat{\mathcal K}}}
\def\kprod{{\mathbb K}}
\def\C{{\mathcal C}}

\def\cprod{{\mathbb C}}
\def\A{{\mathcal A}}
\def\Aphi{\A^{\phi}}
\def\Apsi{\A^{*,\psi}}
\def\Apsit{\A^{*,\psi}_t}
\def\Aopsi{\A^{o,\psi}}
\def\Aopsit{\A^{o,\psi}_t}
\def\m{m{\mathcal F}^+}

\def\meas{\m}
\def\b{b{\mathcal F}^+}

\def\bounded{\b}
\def\L{{{\mathcal L}}}
\def\Linfty{{\L^\infty}}
\def\Linf{\Linfty}
\def\Lw{{{\mathcal L}^1}}
\def\Lo{\L^0}

\def\Lot{\L^0_t}
\def\LoT{\L^0_T}
\def\Lpt{\L^{0,+}_t}
\def\LpT{\L^{0,+}_T}
\def\R{{\mathbb R}}
\def\stratnull{{\mathcal N}}

\def\tK{{-\tilde {\mathcal K}}}

\def\tA{{\tilde\A}}

\def\om{\omega}
\def\Om{\Omega}
\def\Ct{\C_t}

\def\At{\A_t}

\def\D{{\mathcal D}}

\def\E{{\mathbb E}}
\def\P{{\mathbb P}}

\def\Prightarrow{\buildrel \Lo\over\longrightarrow}
\def\Lwto{\buildrel \Lw\over\longrightarrow}
\def\asto{\buildrel \hbox{a.s.}\over\longrightarrow}
\def\defto{\buildrel def\over =}

\def\moo{(M_t)_{t=0,\ldots, T}}

\begin{document}
\bibliographystyle{plain}
\title[NO ARBITRAGE AND CLOSURE RESULTS]{NO ARBITRAGE AND CLOSURE RESULTS \\
FOR TRADING CONES WITH TRANSACTION COSTS}
\author{Saul Jacka}
\address{Department of Statistics, University of Warwick, Coventry CV4
 7AL, UK} \email{s.d.jacka@warwick.ac.uk}
\author{Abdelkarem Berkaoui}
\email{berkaoui@yahoo.fr}
\author{Jon Warren}
\email{j.warren@warwick.ac.uk}

\begin{abstract}
In this paper, we consider trading with proportional transaction
costs as in Schachermayer's paper of 2004. We give a necessary and sufficient
condition for $\A$, the cone of claims attainable from zero
endowment, to be closed. Then we show how to define a revised set of trading prices
in such a way that firstly, the corresponding cone of claims attainable for
zero endowment, $\tA$, does obey the
Fundamental Theorem of Asset Pricing and secondly, if $\tA$ is arbitrage-free then it is the closure
of $\A$. 
We then conclude by showing how to represent
claims.
\end{abstract}

\thanks{The suthors would like to thank two anonymous referees, whose
suggestions and comments have greatly enhanced the exposition and
presentation of this paper.}
\thanks{This research was supported by the grant \lq Distributed Risk Management'
in the Quantitative Finance initiative funded by EPSRC and the
Institute and the Faculty of Actuaries} 
\maketitle \centerline{Revision: {\today}}

\section{Introduction, notation and main results}
\subsection{Introduction}Recollect the Fundamental Theorem of Asset Pricing in finite
discrete time (see, for example, Schachermayer \cite{scha2}): the fact that $\A$, the set of
claims attainable for 0 endowment, is arbitrage-free implies and is
implied by the existence of an Equivalent Martingale Measure; in addition,
$\A$ is closed if it is arbitrage-free.

In \cite{scha1}, Schachermayer showed that the Fundamental Theorem
of Asset Pricing fails in the context of trading with
spreads/transaction costs, by giving an example of an $\A$ which is
arbitrage-free, but whose closure does contain an arbitrage (see also
Kabanov, Rasonyi and Stricker \cite{kabanov2}
and \cite{kabanov3}). Consequently it is of interest to investigate further when the cone $\A$ is
closed, and in cases when it is not, to find descriptions of its closure.

Schachermayer then established (Theorem 1.7 of \cite{scha1}) the
equivalence of two criteria associated with the no-arbitrage
condition for the general set-up for trading with
spreads/transaction costs: that robust no-arbitrage implies and is
implied by the existence of a strictly consistent price process.
Here, robust no-arbitrage means loosely that even with smaller
bid-ask spreads there is no arbitrage, whilst a strictly consistent
price process is one taking values in the relative interior of the
set of consistent prices. In Theorem 2.1 of \cite{scha1} he showed
that the robust no-arbitrage condition implies the closure (in
$\L^0$) of the set of attainable claims.

In this paper we shall first give, in Theorem
\ref{ctheorem}, a simple necessary and sufficient
condition for the set of attainable claims to be closed.
We go on to show, in Theorem \ref{t2}, how to amend the bid-ask spreads so
that the new cone of attainable claims does satisfy the original
Fundamental Theorem (i.e. is either arbitrage-free and closed or admits an
arbitrage). Moreover, we show that in the arbitrage-free case
the new cone is simply the closure of the original cone of attainable
claims. Finally, in
section \ref{penult}, we consider representation of attainable claims and
characterize claims attainable for a given initial endowment.  

\goodbreak
\subsection{Notation and main results}
We are equipped with a filtered probability space
$(\Omega,\F,(\F_t:\;t=0,1,\ldots, T),\P)$. We denote the set of non-negative,
real-valued $\F_t$-measurable random variables by $\meas_t$ and the bounded
non-negative,  real-valued $\F_t$-measurable random variables by
$\bounded_t$. We denote the set of $\R^d$-valued $\F_t$-measurable random variables by
$\Lot$ and the non-negative $\R^d$-valued $\F_t$-measurable
random variables by $\Lpt$. More generally, we denote the set of $\F_t$-measurable random
variables taking values in the (suitably measurable) random set $S$ by $\L^0(S;\F_t)$.

We recall the setup from Schachermayer's paper \cite{scha1} for trading
with $d$ assets. A $d\times d$ matrix, $\Pi$ is said to be a {\em bid-ask
matrix} if
\begin{itemize} 
\item{}$\Pi^{ij}>0$ for all $i,j$; 
\item{}$\Pi^{ii}=1$; 
\end{itemize} 
and
\begin{itemize}
\item{}$\Pi^{ij}\Pi^{jk}\geq \Pi^{ik}$. 
\end{itemize} 

We interpret $\Pi^{ij}$ as the number
of units of asset $i$ required to purchase one unit of asset $j$.

An  adapted $\R^{d\times d}$  process 
$(\bid_t:\;t=0,1,\ldots, T)$ with each  $\bid_t$ being a bid-ask matrix is  
known as a {\em bid-ask process} and gives the time $t$
price for one unit of each asset in terms of each other asset. We assume
that we are given a fixed bid-ask process, $\bid$.

Next we define, for a fixed bid-ask matrix, $\Pi$, the solvency cone,
$\solvD(\Pi)$, as the convex cone in $\R^d$
spanned by the 
canonical basis vectors of $\R^d$, $(e_i)_{1\leq i\leq d}$, together with the
vectors $\Pi^{ij}e_i-e_j$. The solvency cone thus consists of all those holdings which can be
traded to a non-negative holding at the prices specified by $\Pi$.

The cone of portfolios available at price zero under the bid-ask matrix
$\Pi$ is $\tradeD(\Pi)$.

The time $t$ trading cone consists of all those portfolios
(including those attainable by the \lq\lq burning" of assets) which are available at
time $t$ from zero endowment.  A moment's thought will show that the set
of trades which will be available at time $t$ is the convex cone
$\L^0(\tradeD(\bid_t);\F_t)\defto\trade_t$.

The fundamental object of study is the cone of claims attainable from zero endowment, 
which will be   denoted  by $\A$, and is defined to be
$$
(\K_0)+\ldots+(\K_T).
$$
We also consider
$$
\C_t\defto\{X\in \Lot:cX\in \A\,\hbox{ for all }c\in\b_t\}.
$$
We say a few words on the interpretation of $\C_t$ versus $\K_t$.
It is clear that $\trade_t \subseteq \C_t\subseteq\A$, thus
we have the equality
 $$ \A=\C_0+\ldots +\C_T.$$
We can think of $\C_t$ as consisting of those trades which are available
on terms that are known at time $t$ but which may require trading at
later times to be realised.

Although each $\K_t$ is closed in $\Lot$, this is not enough to ensure 
that $\A$ is closed in $\LoT$. In contrast we find 
 the following  necessary and sufficient condition for the closure of $\A$:
\begin{theorem}
\label{ctheorem} $\A$ is closed in $\LoT$ if and only if each $\C_t$
is closed.
\end{theorem}

Let $\bar\A$ denote the closure of $\A$ in $\LoT$. Unlike in a 
classical market, $\A$ can be arbitrage-free, that is to say 
$$
\A \cap \LpT= \{0\},
$$
yet not closed. It is then natural to ask for a description of the closure, $\bar\A$. 
\begin{theorem}
\label{t2}
There is an adjusted bid-ask process $\tbid$ (see Definition \ref{adjusted}) such that the associated cone
of claims $\tA$ satisfies $\A\subseteq\tA\subseteq\bar\A$. Moreover,
either $\tA$ contains an arbitrage or it is
arbitrage-free and closed. In the former case, $\bar A$ also contains an arbitrage, while in the latter case
$$\bar{\A}=\tA.$$
\end{theorem}

\goodbreak
\section{Results on the closedness of $\A$}
As we have remarked already, $\A$ can be arbitrage-free but not closed.
Recall that Schachermayer gives a sufficient condition for the
closedness of $\A$ in terms of robust arbitrage.

Schacheramyer defines the bid-ask spreads as the (random) intervals
$[\frac{1}{\bid^{j,i}_t},\bid^{ij}_t]$, for $i,j\in \{1,\ldots,d\}$ and $t=0,\ldots,T$, and defines
robust no-arbitrage as follows:
\begin{itemize}
\item{}the
bid-ask process $\bid$ satisfies {\it robust no-arbitrage}\ if there is a
bid-ask process $\tilde\bid$ with smaller bid-ask spreads than $\bid$
(i.e. one whose bid-ask spreads almost surely
fall in the relative interiors, in $\R$, of the bid-ask spreads for $\bid$) whose
cone of admissible claims is arbitrage-free.
\end{itemize}

Theorem $2.1$ of Schachermayer \cite{scha1} then states that robust no-arbitrage
implies that the cone $\A$ is closed --- as the remark after the proof
states, the proof relies only on the collection of {\em null strategies} (see Definition \ref{null})
being a closed vector space. However it is easy to find an example
where $\A$ is closed and arbitrage-free but robust no-arbitrage fails.

Consider the following example:
\begin{example}\label{eg1}
Suppose that $T=1$, $d=2$, $\bid_0^{1,2}=1$, $\bid_0^{2,1}=2$ whilst
$\bid^{ij}_1=1$ for each pair $i,j$. Take $\Omega={\mathbb N}$, $\F_0$ trivial
and $\F_1=2^{{\mathbb N}}$ with $\P$ given by $\P(n)=2^{-n}$.

It is immediately clear that robust no-arbitrage cannot hold, since any 
bid-ask process $\tilde\bid$ with smaller bid-ask spreads than $\bid$ must
have $\tilde\bid^{1,2}_0\in(\frac{1}{2},1)$ and $\tilde\bid^{2,1}_1=1$. There is
then an arbitrage in the corresponding cone $\tilde\A$ since $e_2-
\tilde\bid^{1,2}_0e_1+e_1-\tilde\bid^{2,1}_1e_2$ must be  a positive multiple of
$e_1$.
\end{example}

\begin{remark}
With the setup of Example \ref{eg1}, it is clear from the bid-ask prices that
$$
\trade_0\,=\{(x,y):\: x+y\leq 0\hbox{ and }x+2y\leq 0\}
$$
and
$$\trade_1= \{(X,Y)\in\L^0_1:\, X+Y\leq 0\,\, \P\hbox{
a.s.}\}
$$
and so (since $\trade_0\subset\trade_1$ and $\A=\trade_0+\trade_1$)
$$
\A= \{(X,Y)\in\L^0_1:\, X+Y\leq 0\,\, \P\hbox{
a.s.}\}.$$
We can then see that
$\C_0=\{ (x,y):\,
x+y\leq 0\}$,
while
$\C_1=\A=\{(X,Y)\in\L^0_1:\, X+Y\leq 0\,\, \P\hbox{
a.s.}\}$.
\end{remark}
It is tempting to speculate that if $\A$ is not closed, then $\bar\A$
contains an arbitrage. The following example (compare with example 1.3 in
Grigoriev \cite{Grigoriev}) shows that this is false.

\begin{example}
\label{eg3}
Suppose that $T=1$, $d=2$, $\bid_1^{1,2}=1$, $\bid_1^{2,1}=2$ whilst
$\bid^{ij}_0=1$ for each pair $i,j$.
Take $\Omega={\mathbb N}$, $\F_0$ trivial
and $\F_1=2^{{\mathbb N}}$ with $\P$ given by $\P(n)=2^{-n}$.

Then we have 
$$
\bar\A=\{(X,Y)\in 
\L^0_1:\, X+Y\leq 0\,\, \P\hbox{
a.s.}\},
$$ whereas 
$$
\A=\{(X,Y)\in 
\L^0_1:\, X+Y\leq 0\,\, \P\hbox{
a.s.} \hbox{ and }2X+Y\hbox{ is a.s. bounded above}\}.
$$
\end{example}

\begin{lemma}
\label{clemma} For each $t$, $\C_t$ is a convex cone in $\Lot$ and
$$
\A=\C_0+\ldots +\C_T.
$$
\begin{proof}
Convexity for $\C_t$ is inherited from $\A$ as is stability under
multiplication by positive scalars. The decomposition result follows
from the fact that $$\K_t\subseteq \C_t$$ and the fact that
$\C_t\subseteq \A$.
\end{proof}
\end{lemma}
\begin{definition}\label{null}
For any decomposition of $\A$ as a sum of convex cones:
$$
\A=\M_0+\ldots +\M_T\,,$$ we call elements of $\M_0\times\ldots\times
\M_T$ which almost surely sum to 0, null-strategies (with respect to
the decomposition $\M_0+\ldots +\M_T$) and denote the set of them by
$\stratnull(\M_0\times\ldots\times \M_T)$. For convenience we denote
$(\K_0)\times\ldots\times (\K_T)$ by $\kprod$ and
$\C_0\times\ldots\times \C_T$ by $\cprod$.
\end{definition}

In what follows we shall often use the lemma below (Lemma 2 in
Kabanov et al \cite{kabanov3}):
\begin{lemma}
\label{s} Suppose that
$$\A=\M_0+\ldots +\M_T$$
is a decomposition of $\A$ into convex cones with $\M_t\subseteq
\Lot$ and $\b_t\,\M_t\subseteq \M_t$ for each $t$; then $\A$ is closed if
$\stratnull(\M_0\times\ldots\times \M_T)$ is a vector space and each
$\M_t$ is closed.
\end{lemma}

\begin{lemma}\label{decomp2}
Suppose that $\A=\M_0+\ldots +\M_T$, where for each $t$, $\M_t\subseteq
\Lot$ and $\b_t\,\M_t\subseteq \M_t$, then
$$\M_t\subset \C_t.$$

Moreover, for each $0\leq t\leq T$,
\begin{equation}\label{decomp3}
\A_t(\C)\defto \C_0+\ldots +\C_t=\A\cap \Lot.
\end{equation}
\end{lemma}
\begin{proof}
The inclusion $\M_t\subset \C_t$ follows immediately from the fact that
$\M_t\subset \A$; the stability under multiplication by $\b_t$; and the
definition of $\C_t$.

To prove the equality (\ref{decomp3}), suppose $X\in \A\cap\Lot$. Let
$$X=\xi_0+\ldots +\xi_T,$$
be a decomposition of $X$ with $\underline\xi\in\cprod$.
It follows from the fact that $X\in\Lot$ and $\xi_s\in\Lot$ for each
$s<t$ that
$$Y=\xi_t+\ldots +\xi_T\in \Lot.$$

Therefore, it is sufficient to show
that
$$(\C_t+\ldots +\C_T)\cap \Lot\subset \C_t.$$

Now take $Y\in (\C_t+\ldots +\C_T)\cap \Lot$ and $c\in \b_t$: clearly
$cY\in\A\cap\Lot$ and hence, by definition, $Y\in \C_t$.
\end{proof}

We may now give the

\leftline{\em Proof of Theorem \ref{ctheorem}}

First assume that $\A$ is closed and $(X_n)_{n\geq 1}$ is a sequence in
$\C_t$ converging in $\L^0$ to $X$. It follows immediately from the
assumption that $cX_n\Prightarrow cX\in\A$ for all $c\in \b_t$,
hence $X\in\C_t$.

For the reverse implication we shall show that $\stratnull(\cprod)$
is a vector space and the result will then follow from Lemma
\ref{s}.

Now suppose $(\xi_0,\ldots,\xi_T)\in\stratnull(\cprod)$ and $c\in
\b_t$ with almost sure upper bound $B$: then, defining
$$
\zeta_s=B\xi_s
$$
for $s\neq t$ and
$$
\zeta_t=(B-c)\xi_t, $$
it is clear (from the definition of $\C_s$)
that
$$(\zeta_0,\ldots ,\zeta_T)\in \cprod,$$
with
$$\sum_0^T \zeta_s=-c\xi_t.$$
It follows that
$$-c\xi_t\in \A,\, \forall c\in \b_t$$
and so $-\xi_t\in \C_t$ for each $t$ so that $\stratnull(\cprod)$ is
a vector space as required.
\hfill$\square$

\begin{remark}\label{strat}
In the proof above we used the fundamental property of
null strategies: if $(\xi_s)_{0\leq s\leq T}$ is a null strategy
then $-\xi_t\in \C_t$. A null strategy allows one to eliminate friction
in any of its component trades.
In what follows we shall generalize this idea
to more general sequences of strategies.
\end{remark}


\section{A revised fundamental theorem of asset pricing}\label{adjust}
We return to Example \ref{eg3}:

\begin{example}
\label{eg4}
Recall that $T=1$, $d=2$, $\bid_1^{1,2}=1$, $\bid_1^{2,1}=2$ whilst
$\bid^{ij}_0=1$ for each pair $i,j$;
$\Omega={\mathbb N}$, $\F_0$ is trivial
and $\F_1=2^{{\mathbb N}}$ with $\P$ given by $\P(n)=2^{-n}$.

We leave it as an exercise for the reader to show, as claimed above, that 
$\bar \A=\{(X,Y)\in\L^0_1:\, X+Y\leq 0\,\, \P\hbox{
a.s.}\}$ and hence
corresponds to an {\em adjusted bid-ask process}, which is identically
equal to 1. To do so, one may consider the null strategy $\xi$ given by 
$\xi_0=e_1-e_2$ and $\xi_1=e_2-e_1$.

\end{example}

In this section we shall show that $\bar\A$, if arbitrage-free, can
always be represented by some adjusted bid-ask process. However, the next example, which is a minor adaptation
of one of the key examples in Schachermayer \cite{scha1}, shows
that it is necessary to consider more than just null strategies when
seeking the appropriate adjusted prices.

\begin{definition}
We define $\C_t(\bar\A)$ by analogy with $\C_t(\A)$:
$$
\C_t(\bar\A)\defto\{X\in \Lot:cX\in \bar\A\,\hbox{ for all }c\in\b_t\}.
$$
\end{definition}

\begin{example}
Suppose that $T=1$, $d=4$, $\Omega={\mathbb N}$, ${\mathcal F}_0$ is trivial and ${\mathcal F}_1=2^\Omega$.
The bid-ask process at time 0 satisfies $\bid_0^{2,1}=1$, $\bid_0^{4,3}=1$ whilst
$\bid^{ij}_0=3$ for each other pair $i,j$ with $i\neq j$. At time 1,
we have $\bid_1^{1,4}=\omega=\frac{1}{\bid_1^{4,1}}$ and
$\bid_1^{2,3}=\omega=\frac{1}{\bid_1^{3,2}}$, whilst $\bid_1^{4,3}=1$ and $\bid_1^{3,4}=3$. All
other entries are defined implicitly by the criterion
$$\bid^{ij}_1=\min_{i=i_0,\ldots,i_n=j}\bid^{i_0 i_1}_1\ldots \bid^{i_{n-
1} i_n}_1.$$

We shall show that $e_4-e_3,e_2-e_1,e_1-e_2\in \C_1(\bar\A)$ even though there is no {\bf
null} strategy, $\xi$, with $\xi_0=e_1-e_2$ or with $\xi_0=e_2-e_1$ or with $\xi_0=e_3-e_4$.

First, define a sequence of strategies $\xi^N$ as follows:
$\xi^N_0=N(e_1-e_2)$ and
$$
\xi^N_1=\frac{N}{\omega}(e_4-\omega e_1)+(\frac{N}{\omega}-1_{(N\geq
\omega)})(e_3-e_4)+{N}(e_2-\frac{1} {\omega}e_3),
$$
which means that $\xi^N_1=N(e_2-e_1)+1_{(N\geq \omega)}(e_4-e_3)$.

Notice that  $\sum\limits_{t=0}^1 \xi^N_t=1_{(N\geq
\omega)}(e_4-e_3)\Prightarrow e_4-e_3$ as $N\rightarrow\infty$, so we
conclude that $e_4-e_3\in \C_0(\bar\A)$. However, $e_3-e_4\in\trade_1$
and so $((e_4-e_3),(e_3-e_4))$ is null for $\cprod(\bar\A)$ and hence
$e_4-e_3\in\C_1(\bar\A)$.

Now, given an element $X$  of $\bounded_1$ with a.s. bound $B$, consider the strategy
$((N+B)(e_1-e_2)+(e_3-e_4) ,(N+(B-X)(e_2-e_1)+1_{(N+(B-X)\geq \omega)}(e_4-e_3))
$, which sums to $X(e_1-e_2)-1_{(N+(B-X)<\omega)}(e_4-e_3)\Prightarrow
X(e_1-e_2)$ as $N\rightarrow\infty$. This shows that $e_1-e_2\in \C_1(\bar\A)$ and so is also in
$\C_0(\bar\A)$.

Lastly, consider the strategy
$$
(N(e_1-e_2)+(e_3-e_4) ,(N+X))(e_2-e_1)+1_{(N+X\geq
\omega)})(e_4-e_3)),
$$
which sums to
$X(e_2-e_1)-1_{(N+X<\omega)}(e_4-e_3)\Prightarrow X(e_2-e_1)$ as
$N\rightarrow\infty$. This shows that $e_2-e_1\in \C_1(\bar\A)$ and
so is also in $\C_0(\bar\A)$.

It follows that $\bar A$ corresponds to the adjusted bid-ask
process $\tbid$ given, for $t=0$, by:
$\tbid^{1,2}_0=\tbid^{2,1}_0=\tbid^{3,4}_0=\tbid^{4,3}_0=1$,
$\tbid^{i,j}_0=\tbid^{j,i}_0=3$ for $i\in \{1,2\}$ and $j\in \{3,4\}$;
and for $t=1$ by: $\tbid_1^{1,4}=\omega=\frac{1}{\tbid_1^{4,1}}=\tbid_1^{2,3}=\frac{1}{\tbid_1^{3,2}}$,
whilst $\tbid_1^{4,3}=\tbid_1^{3,4}=\tbid_1^{1,2}=\tbid_1^{2,1}=1$. To
see this, notice that the inclusion $\A\subset \tilde\A$ is obvious,
while $\tilde\A$ is closed (by robust no-arbitrage) and the inclusion $\tilde\A\subset \bar\A$
follows from the arguments above.
\end{example}

In order to prove our new version of the Fundamental Theorem
we first define the {\em adjusted} bid-ask process, $\tbid$. This process will
either be equal to the original bid-ask process or frictionless
($\omega$ by $\omega$ and for a given pair $(i,j)$).
\begin{definition}
Given a bid-ask process $\bid$, we define for each $(i,j,t)$\,,
$$
z^{i,j}_t\defto e_j-\bid^{ij}_t\,e_i
$$
and
\begin{equation}\label{rdef}
R^{i,j}_t\defto\{B\in \F_t\;:\;-z^{i,j}_t\,1_B\in \bar{\A}\}\,.
\end{equation}
\end{definition}
\begin{lemma}\label{rdef2}If $B\in\F_t$ then
$$
-z^{i,j}_t\,1_B\in \bar{\A}\Leftrightarrow -z^{i,j}_t\,1_B\in
{\C}_t(\bar{\A}).
$$

\end{lemma}
\begin{proof}Clearly the RHS implies the LHS {\em a fortiori}.

To prove the reverse implication, first note that, by definition of
$\trade_t$,
$$
kz^{i,j}_t\in \K_t \hbox{ for any }k\in
\m_t,
$$
which in turn implies that
\begin{equation}\label{nulla}
kz^{i,j}_t\in \C_t \hbox{ for any }k\in
\m_t,
\end{equation}
since $\trade_t\subset \C_t$.
Now suppose that $c\in \b_t$ with bound $M$, and set
\begin{equation}\label{nullb}
Z\defto c(-z^{i,j}_t1_B)=M(-z^{i,j}_t1_B)+(M-c)z^{i,j}_t1_B.
\end{equation}
The first term on the right hand side of (\ref{nullb}) is in
$\bar\A$ since $M$ is a positive constant, $-z^{i,j}_t1_B$ is in
$\bar\A$ by assumption and $\bar\A$ is a cone. The second term is in
$\bar\A$ by (\ref{nulla}) and, since $\bar\A$ is a convex cone, $Z\in \bar\A$. The
result follows.
\end{proof}

Now observe that the collection $R^{i,j}_t$ is closed under
countable unions. To see this, observe first 
that, since $\bar\A$ is a closed cone,  $R^{i,j}_t$ is closed under
countable, {\em disjoint}, unions. Now notice that, from Lemma
\ref{rdef2}, if $B\in R^{i,j}_t$ and $D\in\F_t$ then $B\cap D\in
R^{i,j}_t$. It follows that if  $(B_n)_{n\geq 1}$ is a sequence in
$R^{i,j}_t$ then $B_n\setminus (\bigcup\limits_{k=1}^{n-1}B_k)=B_n\cap
(\bigcup\limits_{k=1}^{n-1}B_k)^c\in R^{i,j}_t$ and
hence $\bigcup\limits_n B_n\in R^{i,j}_t$. We then deduce, 
by the usual exhaustion argument, that there exists a $\P$-a.s. maximum,
which we denote by $B^{i,j}_t$; that is to say that
$$B\in R^{i,j}_t \hbox{ and } B^{i,j}_t\subseteq B\Rightarrow
\P(B\setminus B^{i,j}_t)=0.
$$
\begin{definition}\label{adjusted}
We define the {\bf adjusted} bid-ask process $\tbid$ as follows :
\[
\hbox{for each pair } i\neq j \hbox{ and for each }t,\;
\tbid^{j,i}_t\defto\frac{1}{\bid^{ij}_t}\,1_{B^{i,j}_t}+\bid^{ji}_t\,1_{(B^{i,j}_t)^c}\,.
\]
\end{definition}

\begin{remark} $\tbid$ need not satisfy the condition:
$$
\tbid^{ik}\leq \tbid^{ij}\tbid^{jk},
$$
but we may still define the corresponding trading cone and apply
Lemma \ref{s}.
\end{remark}

We denote the corresponding trading cones and cone of attainable claims
by $(\tK_t)_{0\leq t\leq T}$ and $\tilde\A$ respectively.
Throughout the rest of the paper we denote $e_j-\tbid^{i,j}_t\,e_i$ by
$\tz^{i,j}_t$.

We now give the

\leftline{\em Proof of Theorem \ref{t2}}
We first show that
$$
\A\subseteq\tA\subseteq\bar{\A},
$$
and then show that $\tA$ is closed if it is arbitrage-free.

{\em Proof that $(\A\subseteq\tA)$:}

Since $\bid^{ij}_t\bid^{ji}_t\geq 1$, it follows from the definition
that $\tbid_t\leq\bid_t$\, for each $t$ and so
$$\K_t\subseteq \tK_t,$$ and hence
$$\A\subseteq \tA.$$

{\em Proof that ($\tA\subseteq\bar\A$):}

we show this by demonstrating that
$$\tK_t\subseteq\bar\A$$
for each $0\leq t\leq T$.

This, in turn, is achieved by showing that
\begin{equation}
\label{new} d\,\tz^{j,i}_t\in \bar\A,\hbox{ for all }d\in\m_t.
\end{equation}
From the definition of the adjusted bid-ask process, we obtain :
\[
\tz^{j,i}_t=-\tbid^{j,i}_t\;z^{i,j}_t\,1_{B^{i,j}_t}
+z^{j,i}_t\,1_{(B^{i,j}_t)^c}\,.
\]
Observe that $ -z^{i,j}_t\,1_{B^{i,j}_t}\in \C_t(\bar{\A})$ by
definition of the set $B^{i,j}_t$ and (\ref{nullb}), so
$$
-d\tbid^{j,i}_t\;z^{i,j}_t\,1_{B^{i,j}_t}\in \C_t(\bar{\A})\subset
\bar{\A},
$$
and
$$
d\,z^{j,i}_t\,1_{(B^{i,j}_t)^c}\in \K_t\subseteq\bar\A$$ by
definition of $\K_t$, so that $d\,\tz^{j,i}_t\in \bar\A$ as
required.

{\em Proof that ($\tA$ is closed if $\tA$ is arbitrage-free):}

We prove this by showing that the nullspace
$\tilde\stratnull\defto\stratnull \left((\tK_0)
\times\ldots\times (\tK_T)\right)$ is a vector space and then appealing to
Lemma \ref{s}.

Let $\xi\in \tilde\stratnull$. Then, defining ${\C}_t(\tilde{\A})$
analogously to ${\C}_t({\A})$, for each $t$ we have, by Remark
\ref{strat}, $-\xi_t\in {\C}_t(\tilde{\A})$, because $\xi$ is null
for $\tilde A$.

Now, since $\xi_t\in\tK_t$ we may write it as
\[
\xi_t=\sum_{i,j}\,\alpha^{i,j}_t\,\tz^{i,j}_t-\sum_k\beta^k_te_k,
\]
for suitable $\alpha^{i,j}_t$ and $\beta^k_t$ in $\b$.
Moreover, $-\xi_t\in {\C}_t(\tilde{\A})$ and since
$\sum_{i,j}\,\alpha^{i,j}_t\,\tz^{i,j}_t\in\tilde\A$ we conclude that
$\sum_k\beta^k_te_k\in\tilde{\A}$. Now, since, by assumption, $\tA$ is arbitrage-free,
we conclude that $\sum_k\beta^k_te_k=0$ a.s., so
\[
\xi_t=\sum_{i,j}\,\alpha^{i,j}_t\,\tz^{i,j}_t,
\]
and consequently $-\sum_{i,j}\,\alpha^{i,j}_t\,\tz^{i,j}_t\in
{\C}_t(\tilde{\A})$. Since ${\C}_t(\tilde{\A})$ is a convex cone and
$\alpha^{i,j}_t\,\tz^{i,j}_t\in \tK_t\subset {\C}_t(\tilde{\A})$ for
each $(i,j)$, we may deduce that, for each pair $(i,j)$:
$$
-\alpha^{j,i}_t\,\tz^{j,i}_t\in {\C}_t(\tilde{\A}).
$$
Now, multiplying by the positive, bounded and $\F_t$-measurable r.v.
$\frac{1}{\alpha^{j,i}_t}1_{\left(\{\alpha^{j,i}_t>\frac{1}{n}\}\cap
(B^{i,j}_t)^c\right)}$, we see that
\[
-z^{j,i}_t\,1_{\left(\{\alpha^{j,i}_t>\frac{1}{n}\}\cap
(B^{i,j}_t)^c\right)}=-\tz^{j,i}_t\,1_{\left(\{\alpha^{j,i}_t>\frac{1}{n}\}\cap
(B^{i,j}_t)^c\right)}\in \tilde{\A}\subset \bar\A.
\]
Then, by definition of the set $B^{j,i}_t$, for each $n$ the subset
$D^{i,j}_t(n)\defto\{\alpha^{j,i}_t>\frac{1}{n}\}\cap
(B^{i,j}_t)^c\subset B^{j,i}_t$. Now, by taking the union over $n$, we see that
$$
D^{i,j}_t\defto\{\alpha^{j,i}_t>0\}\cap
(B^{i,j}_t)^c=\cup_nD^{ij}_t(n)\subset B^{j,i}_t,$$
and we obtain therefore that
$$
\tbid^{j,i}_t=\bid^{ji}_t=\frac{1}{\tbid^{i,j}_t}
$$
on the subset $D^{i,j}_t$. We deduce that
\[
-\tz^{j,i}_t\,1_{D^{i,j}_t}=-z^{j,i}_t\,1_{D^{i,j}_t}=\tbid^{j,i}_t\,\tz^{i,j}_t\,1_{D^{i,j}_t}\in
\tK_t\,,
\]
and
\[
-\tz^{j,i}_t\,1_{\left(\{\alpha^{j,i}_t>0\}\cap
B^{i,j}_t\right)}=\tbid^{j,i}_t\,
z^{i,j}_t\,1_{\left(\{\alpha^{j,i}_t>0\}\cap B^{i,j}_t\right)}\in
-K_t\subset \tK_t\,.
\]
Hence $-\xi_t\in \tK_t$. It follows that $\tilde\stratnull$ is a vector space
as claimed.
\hfill$\square$

%

\section{Decompositions of $\A$, representation and dual
cones}\label{penult}
\subsection{Decompositions of $\A$ and consistent price processes}
We have given a necessary and sufficient condition for $\A$ to be
closed in terms of the ${\C}_t(\A)$ and we have shown how to amend
the bid-ask prices so that the new cone attainable with zero
endowment is $\bar\A$ (if $\bar \A$ is arbitrage-free). It is
natural to ask whether the resulting trading cones $(\tK_t)_{0\leq
t\leq T}$ coincide with the ${\C}_t(\tilde\A)$'s. The following
example shows that this is far from the case.

\begin{example}
Suppose that $T=1$, $d=4$, $\Omega=\{1,2\}$, ${\mathcal F}_0$ is trivial and ${\mathcal F}_1=2^\Omega$.
The bid-ask process at time 0 satisfies  $\bid_0^{4,3}=\bid_0^{4,2}=1$
whilst,
for all other pairs $i\neq j$, $\bid^{ij}_0=4$; the bid-ask process at
time $t=1$ satisfies $\bid_1^{2,1}(1)=4/3=2-\bid_1^{3,1}(1)=2-
\bid_1^{2,1}(2)=\bid_1^{3,1}(2)$ whilst,
for all other pairs $i\neq j$, $\bid^{ij}_1=4$.
By considering the strategy $\xi$ given by $\xi_0=\frac{1}{2}(e_3+e_2)-e_4$ and
$\xi_1=e_1-\frac{1}{2}(e_3+e_2)$, we see that $e_1-e_4\in \A$ and hence
is in ${\C}_0$. Now $\Omega$ is finite so $\A$ is closed and it
is now easy to check that $\tbid=\bid$, yet $e_1-e_4\not\in -K_0$ and so
$\tK_0\neq \C_0$.
\end{example}

In the rest of this section we shall show that nevertheless, the
$\C_t$'s
and their \lq duals' behave like the original trading cones.

Whereas each trading cone, being generated by a finite set of random
vectors, can clearly be identified as $\L^0(S;\F_t)$ for a suitable random cone $S$,
the same is not evidently true of the
$\C_t$s. Thus, we first need some abstract results relating to cones of random
variables.

\begin{remark}
We denote by $\D$, the collection of all closed subsets of $\R^d$.
The standard Borel structure on $\D$, known as the Effros-Borel
structure, and denoted ${\mathcal B}({\mathcal D})$, is as follows: for any set $B$ in $\R^d$ define $\D(B)$ by
$$
\D(B)=\{C\in\D:\,C\cap B\neq \emptyset\},
$$
then ${\mathcal B}(\D)=\sigma(\bid)$, where
$$
\bid=\{\D(B):\, B\hbox{ open in }\R^d\}.
$$
\end{remark}

\begin{definition} Let us consider a map $\Lambda:\Omega\rightarrow \D$. We say that
$\Lambda$ is Effros-Borel measurable if for all open sets
$U\subset\R^d$, we have $\{\omega:\;\Lambda(\omega)\cap U\neq
\emptyset\}\in\F$. We denote by $\Upsilon$, the set of all
Effros-Borel measurable maps. We also  refer to any $\Lam\in \Upsilon$ as a {\em random
closed set}.
\end{definition}

\begin{lemma}
For any $X\in \Lo (\R^d;{\mathcal F})$ and $\Lambda\in \Upsilon$,
\begin{equation}\label{randcon}
(X\in \Lambda)\defto \{\omega:\, X(\omega)\in \Lambda(\omega)\}\in
{\mathcal F}.
\end{equation}
\end{lemma}
\begin{proof}First, by the fundamental measurability theorem of
Himmelberg \cite{Him}, there is a sequence of $\R^d$-valued random
variables $(X_n)_{n\geq 1}$ such that a.s
$$
\Lambda(\omega)=\overline{\{X_n(\omega):\;n\geq 1\}}.
$$
Then, the set $\{\omega:\, X(\omega)\in
\Lambda(\omega)\}=\bigcap\limits_n\bigcup\limits_i\{\omega:\,
|X_i(\omega)- X(\omega)|<\frac{1}{n}\}\in {\mathcal F}$.
\end{proof}
\begin{remark}
In what follows we call a map
$D\in\Upsilon$ with values in the set of closed convex cones in $\R^d$ a {\em random closed cone}.
\end{remark}

\begin{theorem}
\label{abscone} {\bf Abstract closed convex cones theorem.} Let $\C$
be a closed convex cone in $\L^0(\R^d;\F)$, then
\begin{equation}\label{stable}
\C\hbox{ is stable under multiplication by (scalar) elements
of }\b
\end{equation}
iff there is a map $\Lambda\in\Upsilon$ such that
\begin{equation}\label{stable2}
\C= \L^0(\Lambda;\F).
\end{equation}
In this case, the map $\Lam$ is a random closed cone.
\end{theorem}
\begin{proof}The implication (\ref{stable2})$\Rightarrow$ (\ref{stable}) is
obvious.

To prove the
direct implication: we consider the family:
$$
\Upsilon_\C=\{\Gamma\in \Upsilon:\;\L^0(\Gamma;\F)\subset \C\}.
$$

From Valadier \cite{val1} and \cite{val2}, there is an essential
supremum $\Lam\in \Upsilon$ of this family $\Upsilon_\C$, i.e.:
\begin{enumerate}
\item for all $\Gamma\in\Upsilon_{\C}$, we have $\Gamma\subset\Lam$ a.s.;
\item if $\Sigma\in \Upsilon$ is such that for all $\Gamma\in\Upsilon_{\C}$, we have $\Gamma\subset\Sigma$ a.s,
then $\Lam\subset\Sigma$ a.s.
\end{enumerate}
Moreover there is a countable subfamily $(\Gamma^n)_{n\geq
1}\subset\Upsilon_{\C}$ such that $\Lam=\overline{\bigcup_{n\geq
1}\Gamma^n}$ a.s. We want to prove that $\C=\Lo(\Lam;\F)$. To do this, first
we remark that $\C(\Lam)=\overline{\bigcup_{n\geq 1}\C(\Gamma^n)}$.
Then $\Lo(\Lam;\F)\subset \C$ and so $\Lam \in \Upsilon_{\C}$. Now let $\xi\in
\C$ and define the map $\Gamma(\omega)=\Lam(\omega)\cup
\{\xi(\omega)\}$. For $X\in \Gamma$ a.s and $B=\{\xi=X\}$ we have $X
1_{B^c}\in \Lam$ and then $X 1_{B^c} \in \C$ and $X 1_B=\xi 1_B\in
\C$. So $X\in \C$. We deduce that $\Lo(\Gamma;\F)\subset \C$ and then
$\Gamma\in \Upsilon_{\C}$. By the essential supremum property of $\Lam$,
we have $\Gamma\subset \Lam$ and then $\xi \in\Lam$ a.s.

Now suppose that (\ref{stable2}) is satisfied and consider the
sequence $(X_n)_{n\geq 1}$ that generates $\Lam$. For any
$\alpha\in\R^n$, define
$$
Y_{n,\alpha}=\sum_{i=1}^n\,\al_i\,X_i.
$$
Notice that, denoting the non-negative rationals by ${\mathbb Q}_+$, the collection
$$
S\defto{\{Y_{n,\al}:\,\al=(\al_1,\ldots,\al_n)\in {\mathbb Q}_+^n\}}
$$
is countable.

Define the map
$\tLam$ by:
$$
\tLam(\omega)=\overline{\{Y(\om):\;Y\in S\}}^{\R^d}.
$$
From the convex cone property of $\C$, we have each $Y\in \C$
and then, from (\ref{stable2}), $\P(Y\in\Lam)=1$. We deduce that $\tLam\subset
\Lam$ a.s and then (since $X_n\in S$ for each $n$) that $\Lam=\tLam$ a.s. 
\end{proof}
\begin{definition}
Given a closed convex cone $\C$ in $\Lot$ satisfying (\ref{stable})
(with respect to the $\sigma$-algebra $\F$)
we denote
the corresponding random convex cone in (\ref{stable2}) by $\Lambda(\C;\F)$.
\end{definition}
\begin{corollary}
\label{cor1}Suppose that $0\leq p\leq\infty$ and let $\C$ be a convex
cone in $\L^p(\R^d;\F)$ with $\C$ closed in $\L^p(\R^d;\F)$
if $0\leq p<\infty$, and with $\C$
$\sigma(\Linf(\P),\L^1(\P))$-closed if $p=\infty$. Then, $\C$ is
stable under multiplication by (scalar) elements of $\b$
iff there exists a random closed cone $D$ such that
$$
\C=\L^p(D;\F).
$$
\end{corollary}

\begin{proof}First suppose that $0\leq p<\infty$ and
consider $\overline{\C}^{0}\defto \overline{\C}^{\Lo}$,
the closure of $\C$ in $\L^0$.
It is clear that $\overline{\C}^{0}$ inherits stability under multiplication by
$\bounded$ from $\C$ so, by Theorem \ref{abscone},
$$
\overline{\C}^0=\L^0(D;\F),
$$
where $D=\Lambda(\overline{\C}^{0};\F)$.
It suffices then to prove that $\C=\overline{\C}^0\cap \L^p$. The
inclusion $\C\subset\overline{\C}^0\cap \L^p$ is obvious. Now let
$X\in \overline{\C}^0\cap \L^p$, so there exists a sequence $Y^n\in
\C$ which converges a.s to $X$. 
Take a sequence $(\phi_m)_{m\geq 1}$ of continuous functions on $\R$ with
compact support such that $\phi_m$ tends pointwise to 1 as $m\rightarrow \infty$, then,
by the Bounded Convergence Theorem,
$Y^n_m\defto Y^nf_m(|Y^n|)\in \C$ converges to $Y_m\defto
X\,\phi_m(|X|)$ in $\L^p$. So $Y_m\in \C$ and, by letting $m\uparrow
\infty$, we obtain the result that $X\in \C$.

In the case where $p=\infty$, given $X\in
\overline{\C}^{0}\cap\L^\infty$ again take a sequence $(Y^n)$ in $\C$
such that $Y^n\asto X$. Then, for any $f\in \Lw(\R^d;\F)$ and any
$m$, we have that $f.Y^n\phi_m(|Y^n|)\asto f.X\phi_m(|X|)$,
and then $f.Y^n\phi_m(|Y^n|)\Lwto f.X\phi_m(|X|)$ by the
Dominated Convergence Theorem. We conclude that $X\phi_m(|X|)\in
\C$ and hence, again letting $m\uparrow \infty$, we obtain the
inclusion $\overline{\C}^{0}\cap\L^\infty\subset \C$, since $\C$ is
closed in $\sigma(\Linfty,\L^1)$ and hence in $\L^\infty$.
\end{proof}

\begin{lemma}\label{polar}
Let $\C$ be a closed convex cone in $\L^0(\R^d;\F)$, stable under
multiplication by (scalar) elements of $\b$, let $1\leq p<\infty$,  and $\Lambda=\Lambda(\C;\F)$ be as
defined before, then defining
$$
\C^p=\C\cap\L^p,
$$
the polar of $\C^p$ is given by
$$
(\C^p)^*=\L^q(\Lambda^*;\F),
$$
where $q$ is the conjugate of $p$ and $\Lambda^*$ is the polar of $\Lambda$ in $\R^d$.
\end{lemma}

\begin{proof}
This parallels the second half of the proof of Theorem \ref{abscone}.
\end{proof}

\begin{definition}
An adapted sequence of random closed cones in $\R^d$, $(M_t)_{t=0,\ldots, T}$, 
is called a trading decomposition of $\A$ if
$$\A=\Lo(M_0;\F_0)+\ldots +\Lo(M_T;\F_T).$$ 

For such a decomposition,
set $\M_t=\Lo(M_t;\F_t)$ and,
recalling that $\Mprod$ denotes $\M_0\times\ldots\times \M_T$, set
$$
\At(\Mprod)\defto \M_0+\ldots +\M_t.
$$

For any trading decomposition $\moo$, we define a {\em consistent price process}
(with respect to $\moo$) to be
a martingale, $Z$, with $Z_t$ taking values in
$M_t^*\setminus \{0\}$ for each $t$. Thus, a consistent price process is
nothing but a martingale selection of the set-valued process
$(M^*_t\setminus\{0\})$.

\end{definition}

Let $\phi:\Omega\rightarrow (0,1]$ be an $\F_T$-measurable positive random
variable. We denote by $\L^1_{\phi}$ the Lebesgue space associated to
the norm defined by
$$
||f||_{\L^1_{\phi}}\defto\E\{\phi\,|f|_{\R^d}\}\,.
$$
Its dual,
denoted by $\L^{\infty}_{\psi}$, with $\psi=\frac{1}{\phi}$, is associated
with the norm
$$
||f||_{\L^{\infty}_{\psi}}=\hbox{ess sup}\{\psi\,|f|_{\R^d}\}\,.
$$

\begin{theorem}\label{fundament2}
$\bar\A$, the closure of
$\A$ in $\Lo$, is arbitrage-free iff there is a consistent (for some and then for 
any trading decomposition $\moo$ of $\A$) price process
$Z$, and in this case, for every strictly positive $\F_T$-measurable $\phi:\Om\rightarrow
(0,1]$ we may find a consistent price process $Z$ such that $|Z_T|\leq
c\phi$ for some positive constant $c$. In particular, taking $\phi=1$, we
can find a {\em bounded} consistent price process iff $\bar\A$ is closed.
\end{theorem}
\begin{proof}
This follows very closely the proof of Theorem 1.7 (assuming Theorem
2.1) of Schachermayer \cite{scha1}, ignoring references to \lq
robust' and \lq strict'. A sketch proof is as follows: under the
assumption that $\bar\A$ is arbitrage-free, an exhaustion argument
(see \cite{yan}), establishes the existence of a strictly positive
element, $Z$, of the polar to $\bar\A\cap \L^1_\phi$, whilst Lemma
\ref{polar} and the fact that $\M_t\subset\A$ establishes that
$Z_t\defto\E[Z|\F_t]\in \Lambda^*(\M_t;\F_t)$. Conversely, given a
consistent $Z$, we define a frictionless bid-ask process $\hat\bid$
by
$$
\hat\bid^{ij}_t=\frac{Z^j_t}{Z^i_t}.
$$
Taking $Z^1$ as num\'eraire and observing that $\bbQ$ given by
$\frac{d\bbQ}{d\P}$ is an EMM for the corresponding discounted asset
prices, we see, by applying the fundamental theorem for frictionless
trading, that $\hat\A$ is closed and arbitrage-free. Now it is
clear, since $Z$ is a consistent price process, that $\M_t\subset
\hatK_t=\{X\in\Lot:\, Z_t.X\leq 0\hbox{ a.s.}\}$ and hence it
follows that $\bar A$ is arbitrage-free.
\end{proof}

Similar results were proved in Stricker \cite{Stricker}, Jouini and
Kallal \cite{JK}, Schachermayer \cite{scha1} and Grigoriev
\cite{Grigoriev}.

We denote $\A\cap \L^1_\phi$ by $\Aphi$ and by $\Apsi$ its polar cone.
We denote the consistent price processes with $Z_T\in \A^{*,\psi}$ by
$\Aopsi$, and the sets $\{X:\;X=Z_t\hbox{ for some } Z\in \Apsi\}$ and
$\{X:\;X=Z_t\hbox{ for some } Z\in \Aopsi\}$ by $\Apsit$ and $\Aopsit$ respectively.
\begin{remark}
Notice that if $\Aopsi$ is non-empty, then, identifying martingales
with their terminal values, $\Apsi$ is the closure in
$\L^{\infty}_{\psi}$ of $\Aopsi$. This is a standard argument,
following from the fact that if $X\in \Apsi$ and $Y\in \Aopsi$,
then $X+\epsilon Y\in \Aopsi$ for every $\epsilon>0$. It also follows
that $\Apsit$ is the closure in $\L^{\infty}_{\psi}$ of $\Aopsit$.

\end{remark}
\begin{remark}
Note that in Theorem \ref{fundament2}, we do not need to assume that
$\A$ is decomposed as a sum of $\trade_t$'s, but merely that it admits a trading decomposition.
\end{remark}


\goodbreak
\begin{lemma}\label{lem5.3}
Let $X\in \L^1_\phi$. Then the following assertions are equivalent.
\begin{enumerate}
\item $X\in \Ct^\phi\defto\Ct\cap \L^1_{\phi}$.
\item $ X\in \L^1_\phi(\F_t)$ and $Z_t.X\leq 0$ a.s. for all $Z\in \A^o_\psi$.
\item $\E[(W.X)|\,\F_t]\leq 0$ for all $W\in \L^{\infty,+}_\psi$ such that $\E[W|\F_t]\in \A^{0,\psi}_t$.
\end{enumerate}
\end{lemma}
\begin{proof} ($(1)\Rightarrow(2)$)

Clearly, if $X\in\C^\phi_t$, $X\in \L^1_\phi(\F_t)$. Now, for $Z\in
\A^o_\psi$ and $f\in \b_t$ we have:
$$
\E f(Z_t.X)=\E Z_t.(f\,X)=\E Z_T.(f\,X)\leq 0\,,
$$
since $Z_T\in \A^*_\psi$ and $f\,X\in \A^\phi$. Since $f$ is arbitrary it follows that $Z_t.X\leq 0$
a.s.

($(2)\Rightarrow(1)$)

Now let $f\in \b_t$ and $X$ satisfy (2). We need only prove that $fX\in
\A$.

Let $Z\in \A^o_\psi$ then
$$
\E Z_T.(f\,X)=\E Z_t.(f\,X)=\E f(Z_t.X)\leq 0\,.
$$
Therefore, given $Z\in\A^*_\psi$, by taking a sequence $(Z_n)_{n\geq
1}$ in $\A^o_\psi$ converging in $\L^{\infty}_{\psi}$ to $Z$ we
conclude that  $\E Z_T.(f\,X)\leq 0$ and hence $fX\in
\A^\phi\subset\A$.

($(2)\Rightarrow(3)$)

We remark that for $X$
satisfying $(2)$ we have, for every $W\in \L^{\infty,+}_\psi$ such
that $\E[W|\F_t]\in \Aopsit$ and $f\in\b_t$\,,
$$
\E(f\,(W.X))=\E(f\,\E(W|\,\F_t).X)\leq 0.
$$
Since $f$ is an arbitrary element of $\b_t$,
$$
\E[(W.X)|\,\F_t]\leq 0.
$$
($(3)\Rightarrow(2)$)

Take an $X$ satisfying $(3)$. We prove first that $X\in
\L^1_\phi(\F_t)$.

From $(3)$ we deduce that for every $W\in \L^{\infty,+}_\psi$ we
have $\E[(W-\E(W|\,\F_t)).X]=0$ since
$$
\E[(W-\E(W|\,\F_t))|\,\F_t]=0\in \Apsit.
$$
Consequently for every $W\in \L^{\infty,+}_\psi$ we get
$$
\E W.(X-\E(X|\,\F_t))=\E(W-\E(W|\,\F_t)).X=0\,.
$$
Since $W$ is an arbitrary element of $\L^{\infty,+}_\psi$ we may
deduce that $X=\E(X|\,\F_t)$. Let $Z_t\in \Aopsit$\,, then
$$
Z_t.X=\E(Z_t.X|\,\F_t)\leq 0.
$$
\end{proof}

\subsection{Representation}
The following is an easy modification of Theorem 4.1 of
Schachermayer \cite{scha1} and Theorem 4.2 of Delbaen, Kabanov and
Valkeila \cite{DKV}:
\begin{theorem}\label{member}
Suppose that $\theta\in\LoT$ and $\A$ is closed and arbitrage-free. The following are equivalent:
\begin{itemize}
\item[(i)]There is a self-financing process $\eta$ such that
$$\theta\leq \eta_T,$$
i.e. $\theta\in\A$.

\item[(ii)]For every consistent pricing process $Z$ such that the negative
part $(\theta.Z_T)_-$ of the random variable $\theta.Z_T$ is integrable,
we have
$$
\E[\theta.Z_T]\leq 0.
$$
\end{itemize}
\end{theorem}
\begin{proof}
The proof is a much simplified version of the proof of Theorem 4.1 of Schachermayer
\cite{scha1}. We give a sketch of the proof.

{\em (i)$\Rightarrow$(ii)}

It is easy to check that Remark 2.4 of Schachermayer
\cite{scha1} remains valid if we replace the assumption there that
$\bid$ satisfies the robust no-arbitrage assumption by the assumption
that $\A$ is closed and arbitrage-free, or indeed, merely the assumption that there is a
consistent price process. With this change, we have the forward
implication.

{\em (ii)$\Rightarrow$(i)}

Fix $\theta$ and suppose that (i) does not hold. Now choose a $\phi$ such that $\theta\in\L^1_\phi$. 
Note that $\A^\phi$ is a closed, convex cone  in $\L^1_\phi$. Since 
$\theta\not\in \Aphi$, there exists a separating continuous linear
functional $Z\in \L^\infty_\psi$ such that $Z|_{\Aphi}\leq 0$ and
$<Z,\theta>=E[Z.\theta]>0$. It follows from the first of these properties that
$Z_t=E[Z|\F_t]$ is a consistent price process, and then the second shows
that (ii) fails.
\end{proof}

We may now consider representation of elements of $\A$:
\begin{theorem}
Suppose $\theta\in \A^\phi$ and $\eta$ is an adapted $\R^d$-valued
process in $\L^1_\phi$ with $\eta_T=\theta$, and define
$\xi=(\xi_0,...,\xi_T)$ by $\xi_t\defto\eta_t-\eta_{t-1}$ with
$\eta_{-1}\equiv 0$. Then $\xi\in\prod_0^{T}\C_t^\phi$ if and only
if for all $Z\in \A^o_\psi$, the process $M^Z$ defined by
$M^Z_t=\eta_{t-1}.Z_t$\,, is a supermartingale and $M^Z_T\geq
\theta.Z_T$.
\end{theorem}

\begin{proof}Let $\xi\in\prod_0^{T}\C_t^\phi$ and $Z\in
\A^o_\psi$. Then
$$
\E(M^Z_{t+1}|\,\F_t)=\E(\eta_t.Z_{t+1}|\,\F_t)=\eta_t.Z_t=M^Z_t+\xi_t.Z_t\leq
M^Z_t\,,
$$
since $\xi_t\in \C_t^\phi$ and $Z\in\A^o_\psi$. Moreover we
have
$$
M^Z_{T}=\eta_{T-1}.Z_{T}=-\xi_T.Z_T+\theta.Z_T\geq \theta.Z_T\,,
$$
by the same argument. Conversely, we prove that for every $t$\,,
$\xi_t\in \C_t^\phi$: by Lemma \ref{lem5.3} we need to prove that
$Z_t.\xi_t\leq 0$ a.s for every $Z\in \A^o_\psi$ which is the
case since, for $t\leq T-1$,
$$
\xi_t.Z_t=\E(M^Z_{t+1}|\,\F_t)-M^Z_t\leq 0\,,
$$
and for $t=T$ we have
$$
\xi_T.Z_T=\theta.Z_T-M^Z_T\leq 0.
$$
\end{proof}

\begin{conjecture}
We would like to show that
\begin{equation}\label{conj}
\Aphi=\C^\phi_0+\ldots +\C^\phi_T,
\end{equation}
or just that
$$
\Aphi=\overline{\C^\phi_0+\ldots +\C^\phi_T}^{\L^\phi}.
$$
but a proof of either statement eludes us.

We conjecture that (\ref{conj}) is true.
\end{conjecture}

\begin{remark}
We can consider $\eta$'s only defined for $t\leq T-1$ in the theorem
above to obtain the following:
\begin{corollary}
Suppose that $\eta$ is adapted to $(\F_t:0\leq t\leq T-1)$. Then
$\xi\in \prod_0^{T-1}\C^\phi_t$ if and only if the process $M^Z$ is
a supermartingale for all $Z\in D^{0,\psi}$. We may close $\eta$
on the right by $\theta$ if and only if $M^Z_T\geq \theta.Z$ for all
$Z\in D^{0,\psi}$.
\end{corollary}
\end{remark}

\end{document}